# POSTERIOR CONSISTENCY OF GAUSSIAN PROCESS PRIOR FOR NONPARAMETRIC BINARY REGRESSION

BY SUBHASHIS GHOSAL[1] AND ANINDYA ROY

*North Carolina State University and University of Maryland, Baltimore County*

Consider binary observations whose response probability is an unknown smooth function of a set of covariates. Suppose that a prior on the response probability function is induced by a Gaussian process mapped to the unit interval through a link function. In this paper we study consistency of the resulting posterior distribution. If the covariance kernel has derivatives up to a desired order and the bandwidth parameter of the kernel is allowed to take arbitrarily small values, we show that the posterior distribution is consistent in the $L_1$-distance. As an auxiliary result to our proofs, we show that, under certain conditions, a Gaussian process assigns positive probabilities to the uniform neighborhoods of a continuous function. This result may be of independent interest in the literature for small ball probabilities of Gaussian processes.

**1. Introduction.** Consider a binary response variable $Y$ corresponding to a $d$-dimensional covariate $x$. The problem is to estimate the response $p(x) = P(Y = 1|x)$ over the entire covariate space based on an increasing number of observations. We assume that the possible values of the covariate lie in a compact subset $\mathfrak{X} \subset \mathbb{R}^d$. A Bayesian method for estimating $p$ was developed in [4]. A prior on $p$ was induced by the relation $p(x) = H(\eta(x))$, where $\eta$ is a Gaussian process indexed by $\mathfrak{X}$ and $H$ is a known strictly increasing, Lipschitz continuous cumulative distribution function on $\mathbb{R}$. Choudhuri, Ghosal and Roy [4] described algorithms for computing the posterior distribution of $p$ and numerically investigated the properties of the posterior.

In this paper we show consistency of the posterior distribution of $p$, where the prior is assigned through a Gaussian process as in [4]. Statistical procedures are often justified by asymptotics, and posterior consistency plays

Received May 2005; revised November 2005.
[1]Supported in part by NSF Grant DMS-03-49111.

a major role in validating a Bayesian method. The posterior distribution is said to be consistent if the posterior probability of any small neighborhood of the true parameter value converges to one. Because the notion of consistency is dependent on the topology used to define the neighborhoods, one needs to consider an appropriate topology such as the one based on the $L_1$-distance. Because consistency of $p$ is directly related to the distribution of the covariate values, it makes sense to consider $L_1$-distance weighted by the distribution of the covariates or their empirical measure. In the next section we present three different consistency results, for a random covariate with respect to the $L_1$-distance based on the distribution of covariates, for a designed covariate with respect to the $L_1$-distance based on the empirical distribution of the covariate, and finally for a designed one-dimensional covariate with respect to the $L_1$-distance based on Lebesgue measure. The results hold provided that the covariance kernel of the Gaussian process has a certain number of derivatives. We show posterior consistency by verifying prior positivity and entropy (or testing) conditions of the general posterior consistency results of Ghosal, Ghosh and Ramamoorthi [6] or Choudhuri, Ghosal and Roy [3]. An interesting alternative approach to posterior consistency was given in Walker [13].

In the course of our proof we derive two important auxiliary results. First, we show that a Gaussian process assigns positive probability to any uniform neighborhood of a function in the reproducing kernel Hilbert space of the covariance kernel. This result is of significant general interest. Second, we establish a probabilistic bound on the supremum of the derivative of Gaussian processes with covariance kernels that are differentiable up to a certain order.

The complete flexibility in the shape of the sample paths of a Gaussian process makes it an interesting prior for other function estimation problems, such as density estimation or regression function estimation on a bounded interval. The Gaussian process prior was first used in the context of density estimation by Leonard [10] and Lenk [9]. Posterior consistency of the resulting procedure was recently shown by Tokdar and Ghosh [11]. In the context of additive error nonparametric regression, Choi and Schervish [2] established posterior consistency under certain conditions. Following our approach, it seems possible to treat other generalized regression models, such as Poisson regression, in a similar manner, although the test construction method will be problem specific. The natural extension of consistency results will be the characterization of the posterior rate of convergence in the sense of Ghosal, Ghosh and van der Vaart [7]. Some of the results obtained here may be useful for that purpose as well.

The paper is organized as follows. In the next section we state our main results. Positivity of uniform balls under the Gaussian measure is shown in Section 3. In Section 4 we obtain a useful result on the tail of a Gaussian process and its derivatives, which is subsequently used to show that a certain function sieve only spares an exponentially small probability under the Gaussian process prior. Tests with exponentially small error probabilities for testing a function against the complement of an appropriate neighborhood are obtained in Section 5. The results of these sections are used to prove the main theorems in Section 6.

**2. Main results.** In this section we describe the model and the prior and present our main results. Let $Y$ be a binary response corresponding to a $d$-dimensional covariate $x$ and $p(x) = \mathrm{P}(Y = 1|x)$. Let the covariate values belong to a compact subset $\mathfrak{X}$ of $\mathbb{R}^d$. Let $H$ be a known strictly increasing, Lipschitz continuous cumulative distribution function on $\mathbb{R}$ and let $\eta(x) = H^{-1}(p(x))$. A prior on $p(x)$ is induced by a Gaussian process prior on $\{\eta(x) : x \in \mathfrak{X}\}$ with mean function $\mu(x)$ and covariance kernel $\sigma(x, x')$ through the mapping $p(x) = H(\eta(x))$. The covariance kernel is assumed to be of the form

$$(2.1) \qquad \sigma(x, x') = \tau^{-1}\sigma_0(\lambda x, \lambda x'),$$

where $\sigma_0(\cdot, \cdot)$ is a nonsingular covariance kernel and the hyper-parameters $\tau > 0$ and $\lambda > 0$ play the roles of a scaling parameter and (the reciprocal of) a bandwidth parameter, respectively. Let the hyper-priors on $\tau$ and $\lambda$ be $\tau \sim \Pi_\tau$ and $\lambda \sim \Pi_\lambda$, respectively, where $\Pi_\tau$ and $\Pi_\lambda$ are absolutely continuous probability measures on $\mathbb{R}^+$.

Theorem 4 shows that the sample paths of the Gaussian processes can approximate a large class of functions very well and thus, for the purpose of posterior consistency, it is not necessary to consider additional uncertainty in the link function $H$. In fact, the parameter $\tau$ could be taken to be a fixed constant without affecting posterior consistency. However, practical considerations of small sample accuracy suggest putting a suitable prior on $\tau$. Likewise, it is also sensible to consider the possibility of the presence of hyper-parameters in the "trend function" $\mu(x)$; see Remark 2. On the other hand, it is necessary to vary the bandwidth parameter $\lambda$ all over $(0, \infty)$ to obtain posterior consistency.

We shall work with the sieve of response probability functions

$$(2.2) \quad \Theta_n = \Theta_{n,\alpha} = \{p(\cdot) : p(x) = H(\eta(x)), \|D^w \eta\|_\infty < M_n, |w| \leq \alpha\};$$

here and below $D^w \eta$ stands for $(\partial^{|w|}/\partial^{w_1} t_1 \cdots \partial^{w_d} t_d)\eta(t_1, \ldots, t_d)$, $|w| = \sum w_j$, $\alpha$ is some positive integer and $M_n$ is a sequence of real numbers. Let sequences $\lambda_n$ and $\tau_n$ be such that $\Pi_\tau(\tau < \tau_n) = e^{-cn}$ and $\Pi_\lambda(\lambda > \lambda_n) = e^{-cn}$, for some constant $c$. Specific forms of the hyper-priors and the sequences will be discussed later.

Let

(2.3)
$$\mathcal{A} = \left\{ \eta(x) = \sum_{i=1}^{k} a_i \sigma_0(\lambda x, \lambda t_i),\ a_1, \ldots, a_k \in \mathbb{R}, \right.$$
$$\left. t_1, \ldots, t_k \in \mathfrak{X}, k \geq 1, \lambda > 0 \right\}.$$

Then $\bar{\mathcal{A}}$, the closure of $\mathcal{A}$ in the supremum metric, is called the reproducing kernel Hilbert space (RKHS) of $\sigma_0$ (or, equivalently, of $\sigma$). We make the following assumptions.

ASSUMPTION (P). For every fixed $x \in \mathfrak{X}$, the covariance function $\sigma_0(x, \cdot)$ has continuous partial derivatives up to order $2\alpha + 2$, where $\alpha$ is a positive integer to be specified later.

The mean function $\mu(x)$ belongs to the RKHS, $\bar{\mathcal{A}}$, of the covariance kernel $\sigma_0(\cdot, \cdot)$.

The prior $\Pi_\lambda$ for $\lambda$ is fully supported on $(0, \infty)$.

ASSUMPTION (C). The covariate space $\mathfrak{X}$ is a bounded subset of $\mathbb{R}^d$.

ASSUMPTION (T). The transformed true response function $\eta_0$ belongs to $\bar{\mathcal{A}}$.

Assumption (T) implies that $\eta_0$ is uniformly bounded above and below, and hence, $p_0(x) = H(\eta_0(x))$ is bounded away from 0 and 1. In our setup, the free quantities are $\alpha$ and the sequences $M_n, \lambda_n$ and $\tau_n$. We do not require $\Pi_\tau$ and $\Pi_\lambda$ to have specific forms as long as they satisfy some tail conditions specified by the magnitude of the tail cut-off points $\lambda_n$ and $\tau_n$. The quantity $\alpha$ specifies the smoothness of the covariance kernel. The numbers $\alpha$, $M_n$, $\lambda_n$ and $\tau_n$ need to satisfy some interrelation described by the following growth condition.

ASSUMPTION (G). For every $b_1 > 0$ and $b_2 > 0$, there exist sequences $M_n$, $\tau_n$ and $\lambda_n$ such that
$$M_n^2 \tau_n \lambda_n^{-2} \geq b_1 n \quad \text{and} \quad M_n^{d/\alpha} \leq b_2 n.$$

The first part of Assumption (G) will be used to prove exponential decay of the prior probability of the complement of the sieve $\Theta_n$ and the second part will be used to bound the uniform entropy number of $\Theta_n$.

Now we state our main results under different specifications for the covariate values.

2.1. *Random covariate.* Let $P_0^n$ denote the true distribution of the whole data. We first state the posterior consistency result for the case where the covariates arise as a random sample from a distribution $Q$ on $\mathfrak{X}$.

THEOREM 1. *Suppose the random covariate $X$ is sampled from a probability distribution $Q$ on $\mathfrak{X}$. Suppose that Assumptions* (P), (C), (T) *and* (G) *hold. Then for any $\varepsilon > 0$,*

$$\Pi\bigg(p\colon \int |p(x) - p_0(x)|\, dQ(x) > \varepsilon \bigg| Y_1, \ldots, Y_n,\ X_1, \ldots, X_n\bigg) \to 0$$

*in $P_0^n$-probability.*

The covariate measure, $Q$, may be viewed as a fixed quantity or a nuisance parameter. When $Q$ is a fixed quantity, the posterior distribution does not depend on $Q$. Therefore, to evaluate the posterior, we need not actually know $Q$. On the other hand, when $Q$ is treated as an unknown parameter, we need to specify a prior on $Q$. Under the natural assumption that $p$ is unrelated to $Q$, the likelihood for $p$ can be separated out from that of $Q$. Thus, with independent priors, $p$ and $Q$ will be independent a posteriori, and hence, the posterior distribution of $p$ may be obtained without even specifying a prior on $Q$. Note that the posterior for $p$ will be computed the same as in the case of fixed covariates.

If the covariate measure $Q$ permits a Lebesgue density, then we have the following trivial corollary.

COROLLARY 1. *If the covariate distribution $Q$ has Lebesgue density $q$ which is bounded below by some positive constant, then under the conditions of Theorem 1, consistency in the usual $L_1$-distance $\int |p(x) - p_0(x)|\, dx$ holds.*

2.2. *Designed covariate.* In the case of fixed design, often the entire set of covariate values changes with the sample size. This is the case when the covariates are chosen on some equally spaced grids. Thus, the covariate values form a triangular array of the form $\{x_{i,n}, i = 1, \ldots, n\}$, where repetitions are allowed. Let $Q_n$ be the empirical measure of the design points defined as $Q_n = n^{-1} \sum_{i=1}^n \delta_{x_{i,n}}$, where $\delta_x$ denotes the unit mass probability at $x$. Then we have consistency with respect to the $L_1$-distance based on the empirical measure. Such a distance automatically adjusts to the concentration of the covariates and appears to be more intrinsic than the $L_1$-distance with respect to a fixed measure not related to the distribution of the covariate values.

THEOREM 2. *Assume that the covariate values arise from a fixed design. Then under Assumptions* (P), (C), (T) *and* (G), *for any $\varepsilon > 0$,*

$$\Pi\bigg(p\colon \int |p(x) - p_0(x)|\, dQ_n(x) > \varepsilon \bigg| Y_1, \ldots, Y_n\bigg) \to 0$$

in $P_0^n$-probability.

2.3. *One-dimensional covariate.* In the case of a one-dimensional non-random covariate, one can also obtain consistency in the usual $L_1$-sense under an additional assumption on the covariate values. Without loss of generality, assume that the covariate values $x_{i,n}$, $i = 1, \ldots, n$, are in ascending order. Let $S_{i,n} = x_{i+1,n} - x_{i,n}$ be the spacings between consecutive covariate values.

ASSUMPTION (U). Given $\delta > 0$, there exist a constant $K_1$ and an integer $N$ such that, for $n > N$, we have that $\sum_{i: S_{i,n} > K_1 n^{-1}} S_{i,n} \leq \delta$.

The assumption merely states that the measure of the part of the design space where data is sparse is small. Obviously, Assumption (U) is satisfied by any regularly spaced design. If the design was chosen by sampling from a nonsingular distribution, then by the properties of spacings, it can be shown that Assumption (U) holds with probability tending to one.

THEOREM 3. *Suppose that the values of the covariate arise as design points on $\mathfrak{X}$ satisfying Assumption (U) and $\mathfrak{X}$ is a bounded interval of $\mathbb{R}$. Assume that the prior satisfies Assumption (P). The mapping $x \mapsto \eta_0(x)$ and the prior mean $\mu(\cdot)$ are assumed to have two continuous derivatives on $\mathfrak{X}$ and the covariance kernel $\sigma(\cdot, \cdot)$ is assumed to have continuous partial derivatives up to order 6. Assume that $\Pi_\tau$ and $\Pi_\lambda$ are such that $\tau_n^{-1} \lambda_n^4 = O(n)$. Then for any $\varepsilon > 0$,*

$$\Pi\left(p : \int |p(x) - p_0(x)|\, dx > \varepsilon \Big| Y_1, \ldots, Y_n\right) \to 0$$

in $P_0^n$-probability.

2.4. *Remarks.*

REMARK 1. If $\Pi_\tau(\tau < T) = O(e^{-b/T^r})$ as $T \to 0$ and $\Pi_\lambda(\lambda > L) = O(e^{-bL^s})$ as $L \to \infty$, for some $b, r, s > 0$, it then follows that $\tau_n = n^{-1/r}$ and $\lambda_n = n^{1/s}$ satisfy the exponential tail requirement. Then by both assertions of Assumption (G), we have that

$$n \lesssim M_n^2 \tau_n \lambda_n^{-2\alpha} = M_n^2 n^{-(r^{-1} + 2\alpha s^{-1})} \lesssim n^{2\alpha/d} n^{-(r^{-1} + 2\alpha s^{-1})},$$

which implies that $s > d$ and $\alpha \geq (1 + r^{-1})/(2(d^{-1} - s^{-1}))$. In the most favorable case when $r \to \infty$ and $s \to \infty$, we need $\alpha > d/2$. In other words, with the most favorably tailed priors on $\tau$ and $\lambda$, we need to assume the existence of at least $d + 3$ derivatives of the covariance function, with the requirement going up if the tails are thicker. The natural conjugate prior on $\tau$ is a gamma prior, which assigns too much probability to the lower tail and, hence, does not seem to be good in this respect. A better choice would be the inverse gamma prior which corresponds to $r = 1$ and imposes the restriction $\alpha \geq sd/(s-d)$. As there is no natural conjugate prior on $\lambda$, it makes sense to use a prior with a very thin tail, such as $\Pi_\lambda(\lambda > L) \leq e^{-e^L}$. For such thin tails the restriction on $\alpha$ reduces to $\alpha > d$. Then we will need to assume the existence of at least $2d+3$ derivatives of the covariance kernel for our results to hold.

REMARK 2. For practical considerations, it is useful to allow hyperparameters in the mean function $\mu(\cdot)$. For instance, the mean could be taken as a linear combination of a fixed number of functions, so that $\mu(x) = \sum_{j=1}^J \beta_j \psi_j(x)$. A lower degree polynomial is often a good choice. In order to establish posterior consistency under this scenario, one needs to ensure that the tail of the distribution of $\beta$ is thin enough. For instance, if $P(\|\beta\| > B) = O(e^{-c_1 B^{d/\alpha}})$ for some $c_1 > 0$, then $(\|\beta\| > c_2 n^{\alpha/d})$ has exponentially small prior probability for any $c_2 > 0$. Now for any $\beta$ with $\|\beta\| \leq c_2 n^{\alpha/d}$, the complement of the sieve defined by (2.2) continues to have exponentially small prior probability in view of (4.2) below, provided that $c_2$ is chosen small enough depending on $b_2$ in Assumption (G).

REMARK 3. A popular method of prior construction on functions is by expanding the function in a series $\sum \theta_j \psi_j(x)$ and then putting independent $N(0, \tau_j^2)$ priors on the coefficients. Such a prior leads to a Gaussian process prior on the function, where the covariance kernel is $\sigma(x,y) = \sum_j \tau_j^2 \psi_j(x) \psi_j(y)$. Under appropriate differentiability conditions, our results imply posterior consistency at any $p_0 = H(\eta_0)$, where $\eta_0$ belongs to the RKHS of $\sigma$.

**3. Probability of uniform balls.** In this section we establish a property of the support of a Gaussian process which is also of general interest. Let $\{W(t), t \in T\}$ be a Gaussian process indexed by a compact set $T \subset \mathbb{R}^d$, which we take as $[0,1]^d$ without loss of generality. Let the mean function of $W(t)$ be $\mu(t)$ and the covariance kernel be $\sigma(s,t) = \tau^{-1}\sigma_0(\lambda s, \lambda t)$, where $\sigma_0$ is a fixed covariance kernel, and $\tau > 0$ and $\lambda > 0$ are parameters that can possibly vary according to some distribution. We are interested in finding conditions under which $P\{\|W - w\|_\infty < \varepsilon\} > 0$ for some nonrandom function $w(t)$. Such a result was recently obtained by Tokdar and Ghosh [11] using an approach based on conditioning the process at some grid points and then establishing bounds on the conditional mean and variance. Here we provide a shorter proof based on the Karhunen–Loève expansion of the process.

It suffices to show that the result holds for $\lambda$ and $\tau$ varying over a set of positive probability. For our purpose, $\tau$ can be fixed, as the basis function in the Karhunen–Loève expansion is independent of $\tau$. It follows from Lemma 2 of [11] that it suffices to fix $\lambda$ at some suitable value. Thus, with fixed $\tau$ and $\lambda$, $W$ is a Gaussian process. We intend to show, under appropriate conditions on the covariance kernel, that $w$ belongs to the support of the mixture of Gaussian process priors. First, we establish that a Gaussian process, under mild conditions, assigns positive probabilities to the uniform balls around functions in the RKHS of the covariance kernel.

THEOREM 4. *Assume that $\{W(t), t \in T\}$ is a Gaussian process with continuous sample paths having mean function $\mu(t)$ and continuous covariance kernel $\sigma(s,t)$. Assume that $\mu(t)$ and a function $w(t)$ belong to the RKHS of the kernel $\sigma(s,t)$. Then*

$$(3.1) \qquad P\left(\sup_{t \in T} |W(t) - w(t)| < \varepsilon\right) > 0 \qquad \text{for all } \varepsilon > 0.$$

PROOF. We may assume without loss of generality that $\mu(t)$ is the zero function; else we can subtract $\mu(t)$ from $W(t)$ as well as from $w(t)$.

Let $\sum_{k=1}^{\infty} \sqrt{\lambda_i} \xi_i \psi_i(t)$ be the Karhunen–Loève expansion of $W(t)$, so that $\lambda_i$'s are the eigenvalues of the kernel operator $\sigma(s,t)$, $\psi_i(t)$ are the corresponding eigenfunctions and $\xi_i$ are independent $N(0,1)$; see [1], Section III.3. Let $w(t)$ also be represented as $\sum_{k=1}^{\infty} \sqrt{\lambda_i} a_i \psi_i(t)$, where $\sum_{i=1}^{\infty} \lambda_i a_i^2 < \infty$. It follows from Mercer's theorem (cf. Theorem 3.15 of [1]) that the series $\sum_{i=1}^{\infty} \sqrt{\lambda_i} a_i \psi_i(t)$ converges uniformly, and hence, the tail sum is uniformly small.

Bound $|W(t) - w(t)|$ as

$$(3.2) \quad \sup_{t \in T} |W(t) - w(t)| \leq \sup_{t \in T} |W_N(t) - w_N(t)| + \sup_{t \in T} |\bar{w}_N(t)| + \sup_{t \in T} |\bar{W}_N(t)|,$$

where $W_N(t) = \sum_{i=1}^{N} \sqrt{\lambda_i} \xi \psi_i(t)$, $\bar{W}_N(t) = \sum_{i=N+1}^{\infty} \sqrt{\lambda_i} \xi_i \psi_i(t)$, $w_N(t) = \sum_{i=1}^{N} \sqrt{\lambda_i} a_i \psi_i(t)$ and $\bar{w}_N(t) = \sum_{i=N+1}^{\infty} \sqrt{\lambda_i} a_i \psi_i(t)$. Let $\varepsilon > 0$ be given. The second term on the right-hand side of (3.2) is nonrandom and less than $\varepsilon/3$ for $N$ large enough by the uniform convergence.

The basis expansion of the Gaussian process $W_N(t) - w_N(t)$, for any given $N$, has finitely many terms involving i.i.d. $N(0,1)$ variables and continuous function coefficients. Then, the nonsingularity of a normal distribution with nonsingular covariance implies that $P(\sup_{t \in T} |W_N(t) - w_N(t)| < \varepsilon/3) > 0$ for any fixed $N$.

Now if $P(\sup_{t \in T} |\bar{W}_N(t)| < \varepsilon/3) > 0$ for some $N$, then exploiting the independence of $W_N$ and $\bar{W}_N$, it can be easily shown that (3.1) holds. Thus, it suffices to show that $P(\sup_{t \in T} |\bar{W}_N(t)| < \varepsilon) \to 1$ as $N \to \infty$ for any fixed $\varepsilon > 0$. However, as $W(t)$ has continuous sample paths, by assumption, it

follows that

$$P\left(\sup_{t\in T}|\bar{W}_N(t)|\geq \varepsilon\right) = P\left(\sup_{t\in T}|W_N(t)-W(t)|\geq \varepsilon\right)$$
$$\leq \varepsilon^{-2}E\left(\sup_{t\in T}|W_N(t)-W(t)|^2\right),$$

which converges to 0 as $N\to\infty$ by Theorem 3.8 of [1]. This completes the proof. □

Now assume that $\sigma_0$ is bounded away from zero on $T\times T$. Let

(3.3)
$$\mathfrak{C} = \left\{w\in C(T): w(t) = \sum_{i=1}^k a_i\sigma_0(\lambda t, \lambda t_i),\right.$$
$$\left. a_i\in\mathbb{R}, t_i\in T, 1\leq i\leq k, k\geq 1, \lambda > 0\right\}.$$

Let $w_0\in\bar{\mathfrak{C}}$, where $\bar{\mathfrak{C}}$ is the closure of $\mathfrak{C}$. Then by the discussion preceding Theorem 4, it follows that $P(\sup_{t\in T}|W(t)-w_0(t)|<\varepsilon)>0$, where $W$ has the mixture of Gaussian processes distribution discussed there.

For many covariance kernels, $\mathfrak{C}$ is dense in $C(T)$, in which case every continuous function will be in the support of the prior. For example, if $d=1$ and $\sigma_0(s,t)=\psi(s-t)$ for some nonzero, continuous density function $\psi$ on $\mathbb{R}$, then Tokdar and Ghosh [11] showed that $\mathfrak{C}$ is dense in $C(T)$. For higher dimensions, if the covariance kernel is the Kronecker product of one-dimensional kernels in the sense that $\sigma((s_1,\ldots,s_d),(t_1,\ldots,t_d))=\sigma_1(s_1,t_1)\cdots\sigma_d(s_d,t_d)$, where each $\sigma_j$ has RKHS $C(T_j)$, $j=1,\ldots,d$, then Tokdar and Ghosh [11] showed also that the RKHS of $\sigma$ is $C(T_1\times\cdots\times T_d)$. For instance, it follows from there that the kernel $\sigma((s_1,\ldots,s_d),(t_1,\ldots,t_d))=\tau^{-1}\exp[-\sum_{j=1}^d \lambda_j(s_j-t_j)^2]$ on $\mathfrak{X}$, where the $\lambda_j$'s are unrestricted, has RKHS $C(\mathfrak{X})$ for any product type compact domain $\mathfrak{X}$.

## 4. Sieves and tail probabilities.

LEMMA 1. *Let $\Theta_n$ be as defined in (2.2) and Assumptions* (P), (C) *and* (G) *hold. Then $\Pi(\Theta_n^c)\leq Ae^{-cn}$ for some constants $A$ and $c$.*

Because $\Pi_\lambda(\lambda>\lambda_n)\leq e^{-cn}$ and $\Pi_\tau(\tau<\tau_n)\leq e^{-cn}$, it suffices to uniformly bound the probability of $\Theta_n^c$ for given $\lambda\leq\lambda_n$ and $\tau\geq\tau_n$. The lemma will follow from the following result about Gaussian processes which could also be of general interest.

THEOREM 5. *Let $\eta(\cdot)$ be a Gaussian process on $\mathfrak{X}$, a bounded subset of $\mathbb{R}^d$. Assume that the mean function $\mu(\cdot)$ is in $C^\alpha(\mathfrak{X})$ and the covariance kernel $\sigma(\cdot,\cdot)$ has $2\alpha+2$ mixed partial derivatives for some $\alpha \geq 1$. Then $\eta(\cdot)$ has differentiable sample paths with mixed partial derivatives up to order $\alpha$ and the successive derivative processes $D^w\eta(\cdot)$ are also Gaussian with continuous sample paths. Also, the derivative processes are sub-Gaussian with respect to a constant multiple of the Euclidean distance. Further, there exists a constant $d_w$ such that*

$$(4.1) \qquad \mathrm{P}\left(\sup_{x \in \mathfrak{X}} |D^w\eta(x)| > M\right) \leq K(\eta) e^{-d_w M^2/\sigma_w^2(\eta)}$$

*for $w = (w_1, w_2, \ldots, w_d)$, $w_i \in \{0, 1, 2, \ldots, \alpha\}$, $|w| \leq \alpha$ and $\sigma_w^2(\eta) = \sup_{x \in \mathfrak{X}} \mathrm{var}(D^w\eta(x)) < \infty$, $K(\eta)$ is a polynomial in the supremum of the $(2\alpha+2)$-order derivatives of $\sigma$ and the covariance functions of the derivative processes $D^w\eta(x)$ are functions of the derivatives of the covariance kernel $\sigma(\cdot,\cdot)$.*

PROOF. We may assume, without loss of generality, that the mean function is identically zero, because for $M$ sufficiently large,

$$(4.2) \qquad \begin{aligned} \mathrm{P}(\eta : \|\eta\|_\infty > M) &\leq \mathrm{P}(\eta : \|\eta - \mu\|_\infty > M - \|\mu\|_\infty) \\ &\leq \mathrm{P}(\eta : \|\eta - \mu\|_\infty > M/2). \end{aligned}$$

First we show that the process constructed by taking the partial derivative of $\eta$ with respect to the $j$th component, $D_j\eta(\cdot)$, is again a Gaussian process with continuous sample paths and covariance kernel $D_j^2\sigma(\cdot,\cdot)$. Here and below, $D_j^2$ is the partial derivative operator with respect to the $j$th components of both arguments of $\sigma$, that is, $D_j^2\sigma((s_1,\ldots,s_d),(t_1,\ldots,t_d)) = (\partial^2/\partial s_j\, \partial t_j)\sigma((s_1,\ldots,s_d),(t_1,\ldots,t_d))$. According to our general notation, $D_j\eta(\cdot) = D^w\eta(\cdot)$ and $D_j^2\sigma(\cdot,\cdot) = D^wD^w\sigma(\cdot,\cdot)$ componentwise, where $w = e_j$, the $d$-dimensional vector with one at the $j$th place and zeros elsewhere.

To show $D_j\eta$ is again a Gaussian process, we need to investigate the path properties of the one-parameter Gaussian process obtained by letting the $j$th component vary and holding all other $d-1$ parameters fixed. For notational simplicity, we suppress the dependence of the process on the other $d-1$ parameters. By Section 9.4 of [5] (a version of) $\eta(\cdot)$ has continuously differentiable sample paths if

$$|\Delta_h \Delta_h D_j^2\sigma(s,t)| \leq \frac{C}{|\log|h||^a} \qquad \text{as } h \to 0$$

for some $C > 0$ and $a > 3$, where

$$\begin{aligned} \Delta_h \Delta_h D_j^2\sigma(s,t) &= D_j^2\sigma(s+he_j, t+he_j) - D_j^2\sigma(s+he_j, t) \\ &\quad - D_j^2\sigma(s, t+he_j) + D_j^2\sigma(s,t). \end{aligned}$$

Because $\sigma(\cdot,\cdot)$ has bounded mixed partial derivatives of at least up to fourth order, the above condition is trivially satisfied and the process $\eta(\cdot)$, as a process with respect to the $j$th coordinate, has continuously differentiable sample paths. The limit of a sequence of multivariate normal variables is again a multivariate normal, and $D_j \eta(t) = \lim_{h \to 0} (\eta(t + he_j) - \eta(t))/h$.

It follows that $D_j \eta(\cdot)$ is a Gaussian process. Moreover,

$$\mathrm{E}(D_j\eta(s) - D_j\eta(t))^2 = \lim_{h \to 0} \mathrm{E}\left\{\frac{\eta(s+he_j) - \eta(s) - \eta(t+he_j) + \eta(t)}{h}\right\}^2.$$

This follows by the uniform integrability of $(\eta(t+he_j) - \eta(t))^2/h^2$, which is a consequence of the fact

$$\mathrm{E}\left(\frac{|\eta(t+he_j) - \eta(t)|}{h}\right)^4$$
$$= 3\left(\frac{|\sigma(t+he_j, t+he_j) + \sigma(t,t) - 2\sigma(t+he_j, t)|}{h^2}\right)^2 \leq 3B_0^2 < \infty,$$

for some constant $B_0$. Then, the intrinsic semimetric for the partial derivative process is given by

$$\mathrm{E}(D_j\eta(s) - D_j\eta(t))^2 = \lim_{h \to 0} \mathrm{E}\{|\eta(s+he_j) - \eta(s) - \eta(t+he_j) + \eta(t)|\}^2/h^2$$
$$= \lim_{h \to 0} h^{-2}\{\sigma(s+he_j, s+he_j) + \sigma(t+he_j, t+he_j)$$
$$+ \sigma(s,s) + \sigma(t,t) + 2\sigma(s+he_j, t)$$
$$+ 2\sigma(s, t+he_j)$$
$$- 2(\sigma(s+he_j, t+he_j) + \sigma(t, t+he_j)$$
$$+ \sigma(s,t) + \sigma(s, s+he_j))\}.$$

Using the symmetry of the covariance function and by Taylor's expansion, we have, after simplification, that

$$\begin{aligned}(4.3) \quad \mathrm{E}(D_j\eta(s) - D_j\eta(t))^2 &= D_j^2\sigma(s,s) + D_j^2\sigma(t,t) - 2D_j^2\sigma(s,t) \\ &= \sigma^*(s,t), \quad \text{say.}\end{aligned}$$

As the covariance function has bounded mixed partial derivatives up to order $2\alpha + 2$, by Taylor's expansion, we have for $C = \sup\{|D_j^2 \sigma^*(s,t)| : s,t\}$ that

$$(4.4) \quad \mathrm{E}(D_j\eta(s) - D_j\eta(t))^2 \leq C\|s - t\|^2.$$

Thus, the partial derivative process is sub-Gaussian with respect to a constant multiple of the Euclidean distance. Note that from (4.4), the covariance kernel for $D_j\eta(\cdot)$ is given by

$$(4.5) \quad \mathrm{cov}(D_j\eta(s), D_j\eta(t)) = D_j^2\sigma(s,t).$$

Further, as the kernel $D_j^2\sigma(s,t)$ has at least two mixed derivatives with respect to each component, it follows from [5] that the multi-indexed process $D_j\eta$ has continuous sample paths. Thus, the sample paths of $\eta$ are, with probability one, continuously differentiable with respect to each argument.

Replacing $\eta$ by $D_j\eta$, the mixed partial derivative process $D_kD_j\eta$ is again a Gaussian process and is sub-Gaussian with respect to a constant multiple of the Euclidean distance. In general, the derivative process $D^w\eta(s)$ for $|w| \leq \alpha$ is a Gaussian process with covariance kernel $D^wD^w\sigma(\cdot,\cdot)$.

Now $\sigma_w^2(\eta) = \sup\{\text{var}(D^w\eta(s)) : s \in \mathfrak{X}\} < \infty$. We thus have $N(\varepsilon, \mathfrak{X}, \|\cdot\|) \leq C'\varepsilon^{-d}$ and $N(\varepsilon, \mathfrak{X}, \|\cdot\|_\rho) \leq C''\varepsilon^{-d}$ for some constants $C'$ and $C''$ depending on the measure of the set $\mathfrak{X}$ and the kernel $\sigma$. Here $N$ stands for the covering number, $\|\cdot\|$ is the Euclidean distance and $\|\cdot\|_\rho$ is the intrinsic semi-metric of the derivative process $D^w\eta(s)$. The result then follows by applying Proposition A.2.7 of [12], page 442, to each of the derivative processes. $\square$

To complete the proof of Lemma 1, consider the kernel of the form $\sigma(s,t) = \tau^{-1}\sigma_0(\lambda s, \lambda t)$. Let $\xi$ be a process with a fixed covariance kernel $\sigma_0$. Then the mixed derivative processes $D^w\xi(\cdot)$ up to order $\alpha$ have uniformly bounded variances. Now for $\lambda \leq \lambda_n$, $\tau \geq \tau_n$ and $|w| = \alpha$,

$$\sigma_w^2(\eta) = \tau^{-1}\lambda^{2\alpha}\sigma_w^2(\xi) \leq \tau_n^{-1}\lambda_n^{2\alpha}\sigma_w^2(\xi).$$

The rest of the proof now follows easily from Theorem 5 because the contribution from $K(\eta)$ grows only polynomially in $n$.

## 5. Entropy bounds and test construction.

LEMMA 2. *The $\varepsilon$-covering number $N(\varepsilon, \Theta_n, \|\cdot\|_\infty)$, in the supremum norm, of $\Theta_n$ defined by (2.2), is given by $\log N(\varepsilon, \Theta_n, \|\cdot\|_\infty) \leq KM_n^{d/\alpha}\varepsilon^{-d/\alpha}$ for some constant $K$.*

PROOF. The result follows immediately from Theorem 2.7.1 of [12], page 155. $\square$

LEMMA 3. *Let $\nu$ be a finite measure on $\mathfrak{X}$ and let $\psi_1$ and $\psi_2$ be measurable functions such that $0 \leq \psi_1, \psi_2 \leq M$ and $\int |\psi_1 - \psi_2| \, d\nu > (1+\nu(\mathfrak{X}))\varepsilon$ for some $M$, $\varepsilon > 0$. Then $\nu\{x : |\psi_1(x) - \psi_2(x)| > \varepsilon\} \geq \varepsilon/M$.*

PROOF. By the given condition,

$$(5.1) \quad \begin{aligned}(\nu(\mathfrak{X}) + 1)\varepsilon &\leq \int_{x:|\psi_1(x)-\psi_2(x)|>\varepsilon} |\psi_1(x) - \psi_2(x)| \, d\nu(x) \\ &\quad + \int_{x:|\psi_1(x)-\psi_2(x)|\leq\varepsilon} |\psi_1(x) - \psi_2(x)| \, d\nu(x) \\ &\leq M\nu\{x:|\psi_1(x)-\psi_2(x)|>\varepsilon\} + \varepsilon\nu(\mathfrak{X}).\end{aligned}$$

The result now follows by rearranging the terms. □

Applying Lemma 3 to $\psi_1 = p, \psi_2 = p_0$ and $\nu = Q_n$, we obtain

$$(5.2) \qquad I_{n,p} = \#\{x_i : |p(x_i) - p_0(x_i)| > \varepsilon\} \geq K'n$$

for some $K'$. Let $A_p^+ = \{x : p(x) > p_0(x) + \varepsilon\}$ and $A_p^- = \{x : p(x) < p_0(x) - \varepsilon\}$. Then either $A_p^+$ or $A_p^-$ contains at least $K'n/2$ points. For definiteness, assume that $m = m_n = \#I_{n,p}^+ \geq K'n/2$, where $I_{n,p}^+ = \{i : x_i \in A_p^+\}$. For a given $p$, to test the simple null $p_0$ against the simple alternative $p$, we construct a test based on the observations corresponding to only those design points which are in $A_p^+$. Then (5.2) asserts that there is no loss of order of the number of indices. The next lemma is stated in a general framework and shows how to construct such a test.

LEMMA 4. *Let $Y_j$ be independent Bernoulli variables with $\mathrm{P}(Y_j = 0) = \mu_j$, $j = 1, \ldots, m$. Consider testing $H_0 : \mu_j = \mu_{0j}$ against $H_1 : \mu_j = \mu_{1j}$, where $\mu_{1j} > \mu_{0j} + \varepsilon$ for all $j$ and $0 < \varepsilon_0 < \mu_{0j} < 1 - \varepsilon_0 < 1$; here $\varepsilon_0 > 0$ and $\varepsilon > 0$ do not depend on $m$ and $\varepsilon < \varepsilon_0$. Consider the test $\Psi_m = \mathbb{1}\{\sum_{j=1}^m (Y_j - \mu_{0j}) > m\varepsilon/2\}$. Then for all sufficiently large $m$,*

$$(5.3) \qquad \mathrm{E}_{P_0}(\Psi_m) \leq e^{-m\varepsilon^2/2}, \qquad \mathrm{E}_{P_1}(1 - \Psi_m) \leq e^{-m\varepsilon^2/2},$$

*where $P_0$ and $P_1$ are respectively the probability measures under the null and the alternative.*

REMARK 4. The above lemma also holds if $\mu_{1j} < \mu_{0j} - \varepsilon$ for all $j = 1, \ldots, m$ if the test $\Psi_m$ is defined as one that rejects $H_0$ for $\sum_{j=1}^m (Y_j - \mu_{0j}) < -m\varepsilon/2$.

PROOF OF LEMMA 4. By Hoeffding's inequality ([8], Theorem 1),

$$\mathrm{E}_{P_0}(\Psi_m) = P_0(\bar{Y}_m - \mathrm{E}_{P_0}\bar{Y}_m > \varepsilon/2) \leq e^{-m\varepsilon^2/2}.$$

For all sufficiently large $m$, another application of Hoeffding's inequality gives

$$\begin{aligned}\mathrm{E}_{P_1}(1 - \Psi_m) &= P_1(\bar{Y}_m - \mathrm{E}_{P_0}\bar{Y}_m \leq \varepsilon/2) \\ &= P_1\left((\bar{Y}_m - \mathrm{E}_{P_1}\bar{Y}_m) + m^{-1}\sum_{j=1}^m (\mu_{1j} - \mu_{0j}) \leq \varepsilon/2\right) \\ &\leq P_1((\bar{Y}_m - \mathrm{E}_{P_1}\bar{Y}_m) \leq -\varepsilon/2) \leq e^{-m\varepsilon^2/2}.\end{aligned}$$

This completes the proof. □

For a given $p$, we consider the test $\Psi_{n,p}$ which rejects the simple null $p_0$ against the simple alternative $p$ if

$$(5.4) \qquad m^{-1} \sum_{i \in I_{n,p}^+} (Y_i - p_0(x_i)) > m\varepsilon/2.$$

By Lemma 4 above, the test satisfies (5.3) for a simple alternative $p$.

To remove the dependence on $p$, we use the standard technique of covering a set by small balls and estimating the covering numbers. Note that, for a fixed $\varepsilon > 0$, if $i \in A_p^+$ and $\|p^* - p\| < \varepsilon/2$, then

$$(5.5) \qquad p^*(x_i) - p_0(x_i) \geq p(x_i) - p_0(x_i) - \|p - p^*\|_\infty > \varepsilon/2.$$

Therefore, applying Lemma 4 (with $\varepsilon$ replaced by $\varepsilon/2$), we obtain a test $\Psi_{n,p}$ such that $\mathrm{E}_{p_0} \Psi_{n,p} \leq e^{-m\varepsilon^2/8}$ and $\mathrm{E}_{p^*}(1 - \Psi_{n,p}) \leq e^{-m\varepsilon^2/8}$.

With $N = N(\varepsilon/2, \Theta_n, \|\cdot\|_\infty)$, get $p_1, \ldots, p_N \in \Theta_n$ with the property that, for any $p$, there exists a $p_j \in \Theta_n$ such that $\|p - p_j\|_\infty < \varepsilon/2$. Consider the test $\Phi_n = \max(\Psi_{n,p_j}, j = 1, \ldots, N)$. Then

$$(5.6) \qquad \mathrm{E}_{p_0} \Phi_n \leq \sum_{j=1}^N \mathrm{E}_{p_0} \Psi_{p_j,n} \leq N e^{-m\varepsilon^2/8} = \exp(\log N - m\varepsilon^2/8).$$

If $p \in \Theta_n$, choose $j$ such that $\|p - p_j\|_\infty < \varepsilon/2$. Then

$$(5.7) \qquad \mathrm{E}_p(1 - \Phi_n) \leq \mathrm{E}_p(1 - \Psi_{p_j,n}) \leq e^{-m\varepsilon^2/8}.$$

Then by Lemma 2, for any given constant $b_2' > 0$ we can choose a sufficiently small $b_2$ and $M_n$ satisfying Assumption (G) with $b_2$, such that $\log N \leq b_2' n$. Then choosing $m = m_n$ of order $n$, the test $\Phi_n$ satisfies the requirement

$$(5.8) \qquad E_{p_0} \Phi_n \leq e^{-c'n}, \qquad E_p(1 - \Phi_n) \leq e^{-c'n}$$

for some constant $c'$.

## 6. Proof of the main theorems.
Now we prove Theorems 1–3.

PROOF OF THEOREM 1. We are considering the model

$$(6.1) \quad Y_i | X_i \stackrel{\mathrm{ind}}{\sim} \mathrm{Binomial}(1, p(X_i)), \qquad X_i \stackrel{\mathrm{ind}}{\sim} Q, \qquad i = 1, 2, \ldots, n.$$

Then the joint density of $X$ and $Y$ with respect to the product of $Q$ and the counting measure on $\{0, 1\}$, say, $\aleph$, is given by $f(x, y) = p(x)^y (1 - p(x))^{1-y}$. The corresponding true joint density is $f_0(x, y) = p_0(x)^y (1 - p_0(x))^{1-y}$. Recall that by Assumption (T), $\varepsilon_0 < p_0(x) < 1 - \varepsilon_0$ for some $\varepsilon_0 < 1/2$. This implies that $f_0(x, y) > \varepsilon_0$. Also observe that $\int |f_1 - f_2| \, d\aleph \, dQ = 2 \int |p_1 - p_2| \, dQ$ and $\int f_0 \log(f_0/f) \, d\aleph \, dQ = \int p_0 \log(p_0/p) \, dQ + \int (1 - p_0) \log((1 - p_0)/(1 - p)) \, dQ$.

We verify the conditions given in Theorem 2 in [6]. It may be noted that although their result is stated for Lebesgue densities on $\mathbb{R}$, it is valid for densities in any measure space.

We first show that $\Pi\{f: \int f_0 \log(f_0/f) < \varepsilon\} > 0$ for all $\varepsilon > 0$, where $\Pi$ is the prior for $f$, or equivalently,

$$\Pi\bigg\{p: \int p_0 \log \frac{p_0}{p} \, dQ + \int (1-p_0) \log \frac{1-p_0}{1-p} \, dQ < \varepsilon\bigg\} > 0 \qquad \text{for all } \varepsilon > 0.$$

We shall use the following lemma which follows easily from Taylor's expansion.

LEMMA 5. *Let $0 < \varepsilon_0 < \frac{1}{2}$ and $\varepsilon_0 < \alpha, \beta < 1 - \varepsilon_0$. Then there exists a constant $L$ depending only on $\varepsilon_0$ such that*

$$\alpha \bigg(\log \frac{\alpha}{\beta}\bigg)^m + (1-\alpha)\bigg(\log \frac{1-\alpha}{1-\beta}\bigg)^m \leq L(\alpha - \beta)^2, \qquad m = 1, 2.$$

Let $B = \{p: \|p - p_0\|_\infty < \frac{1}{2}c\}$, where $c = \inf\{\min(p_0(x), 1 - p_0(x)): 0 \leq x \leq 1\} > 0$ and $\|p - p_0\|_\infty = \sup\{|p(x) - p_0(x)|: 0 \leq x \leq 1\}$. If $p \in B$, then it follows from Lemma 5 that

(6.2) $\quad \max\bigg(\int p_0 \log \frac{p_0}{p} \, dQ, \int (1-p_0) \log \frac{1-p_0}{1-p} \, dQ\bigg) \leq L\|p - p_0\|_\infty^2.$

Hence, it suffices to show that $\Pi(p: \|p - p_0\|_\infty < \varepsilon) > 0$ for every $\varepsilon > 0$. Because $p_0(\cdot) = H(\eta_0(\cdot))$ and the function $u \mapsto H(u)$ is bounded and Lipschitz continuous, it is enough to show that

(6.3) $\qquad\qquad \Pi(\eta: \|\eta - \eta_0\|_\infty < \varepsilon) > 0 \qquad \text{for every } \varepsilon > 0.$

The result now follows from Theorem 4.

To verify the entropy condition of Theorem 2 in [6], let $\beta > 0$ be given. We consider the sieve $\mathcal{F}_n = \{f(x, y) = p(x)^y (1 - p(x))^{1-y}: p \in \Theta_n\}$, where $\Theta_n$ is defined in (2.2) with $M_n = bn^{\alpha/d}$ and $b > 0$ is a constant to be chosen sufficiently small. By Lemma 2, for some constant $K$, we have that $\log N(\varepsilon, \Theta_n, \|\cdot\|_\infty) \leq K\varepsilon^{-d/\alpha} b^{d/\alpha} n$. Now choosing $b < (\beta/K)^{\alpha/d}\varepsilon$, we can ensure that $\log N(\varepsilon, \Theta_n, \|\cdot\|_\infty) < n\beta$.

Finally, from Lemma 1, we have that $\Pi(\Theta_n^c)$ is exponentially small. This completes the proof. $\square$

PROOF OF THEOREM 2. For nonrandom covariates, the observations are independent, nonidentically distributed. We shall apply Theorem 2 of [3]. The prior positivity condition follows essentially by the same arguments used in the random covariate case. Consider the sieve $\Theta_n$ defined by (2.2). The condition (A3)(iii) of Theorem 2 of [3] holds by Lemma 1. To verify their conditions (A3)(i) and (A3)(ii), we need to show that there exist exponentially consistent tests for testing $H_0 : p = p_0$ against an alternative $H_a : \|p - p_0\|_{1,Q_n} > \varepsilon$ for all $\varepsilon > 0$. The test constructed in (5.8) satisfies the required conditions. □

PROOF OF THEOREM 3. In this case we verify the conditions of Theorem 2 in [3] for the sieve defined by (2.2) with $\alpha = 2$. In view of the proof of our Theorem 2, the only condition that needs to be additionally verified is the testing condition for the usual $L_1$-distance. To construct the required sequence of tests, we estimate the number of covariate values where the true probability function $p_0$ and the alternative $p$ differ by at least a specified amount. The following lemma, where we assume without loss of generality that $\mathfrak{X} = [0, 1]$, estimates that number. As the number is at least a fraction of $n$, it follows that required tests can be constructed as in the proof of Theorem 2. □

LEMMA 6. *For any $p \in \Theta_{n,2}$ such that $\int |p(x) - p_0(x)| \, dx > 5\varepsilon$, (5.2) holds.*

PROOF. For a given function $h$, let $b(x, k, h) = k \sum_{j=1}^{k} (\int_{(j-1)/k}^{j/k} h(t) \, dt) \times \binom{k-1}{j-1} x^{j-1}(1-x)^{k-j}$ stand for the corresponding Bernstein polynomial of order $k$. Then it is well known (and easy to see) that $\sup\{|h(x) - b(x, k, h)| : 0 \leq x \leq 1\} \leq Ak^{-1} \sup\{|h''(x)| : 0 \leq x \leq 1\}$, where $A$ is an absolute constant. Under given assumptions, the choice $M_n \sim \gamma n^{1/2}$ satisfies Assumption (G) for a sufficiently small $\gamma$. Let $\gamma'$ be a given small constant and let $\gamma$ be sufficiently small such that $\gamma < \gamma' \varepsilon / (2A)$. Choose a sequence $k_n \sim \gamma' n$. Let $b(x) = b(x, k_n, p)$ and $b_0(x) = b(x, k_n, p_0)$. Therefore, it follows from the above discussion that, for any $p \in \Theta_{n,2}$, we have that $\|p - b\|_\infty \leq A(\gamma' n)^{-1} \gamma n < \varepsilon/2$, and $\|p_0 - b_0\|_\infty \leq A(\gamma' n)^{-1} \gamma n < \varepsilon/2$. Hence,

$$|p(x) - p_0(x)| > |b(x, k_n, p) - b(x, k_n, p_0)| - \|p - b\|_\infty - \|p_0 - b_0\|_\infty > \varepsilon,$$

if $x \in B_p = \{t : |b(t, k_n, p) - b(t, k_n, p_0)| > 2\varepsilon\}$. Therefore, $B_p \subset A_p := \{x : |p(x) - p_0(x)| > \varepsilon\}$, and hence, it suffices to show that the assertion with $I_{n,p}$ replaced by $I'_{n,p} = \#\{i : x_i \in B_p\}$.

Clearly, $0 \leq b(x), b_0(x) \leq 1$. Also, $\|b - b_0\|_1 \geq \|p - p_0\|_1 - \|p - b\|_\infty - \|p_0 - b_0\|_\infty > 4\varepsilon$. Applying Lemma 3 to the pair $b$ and $b_0$ and $\nu$, the Lebesgue measure on $[0, 1]$, and $\varepsilon$ replaced by $2\varepsilon$, we obtain that $|B_p| > 2\varepsilon$.

Now, as $b$ and $b_0$ are polynomials of order at most $k_n$, the set $B_p$ is at most a union of $k_n$ intervals. Let $J_1, J_2, \ldots$ be these intervals. Find $K_1$ such that Assumption (U) holds for $\delta = \varepsilon$. Call a spacing interval $(x_{i,n}, x_{i+1,n})$ of type I if $S_{i,n} \leq K_1/n$. For any value of $j$, let $J_j^*$ be the union of type I spacing intervals $(x_{i,n}, x_{i+1,n})$ that are completely contained in $J_j$. Note that at most two type I spacing intervals may be partially contained in $J_j$ for any $j$, which has total length bounded by $2K_1/n$. Put $B_p^* = \bigcup_j J_j^*$. Thus by Assumption (U), $|B_p \cap (B_p^*)^c| < \varepsilon + 2K_1 k_n/n$ and hence $|B_p^*| > \varepsilon - 2K_1 k_n/n$. For $j = 1, 2, \ldots$, let $R_j$ be the number of type I spacing intervals completely contained in $J_j$. Considering the possibility that $J_j^*$ may not contain its end points, we find that $J_j^*$ contains at least $R_j - 2$ design points, and hence $B_p^*$ contains at least $\sum_j R_j - 2k_n$ design points. To estimate $\sum_j R_j$, note that $\varepsilon - 2K_1 k_n/n \leq \lambda(B_p^*) \leq \sum_j R_j K_1/n$, and hence $\sum_j R_j \geq n\varepsilon/K_1 - 2k_n$. Hence $B_p$ contains at least $(n\varepsilon/K_1) - 4k_n$ points, which is greater than $n\varepsilon/(2K_1)$ if we choose $\gamma' < \varepsilon/(8K_1)$. $\square$

**Acknowledgment.** Suggestions of the referees have improved the exposition of the paper.

## REFERENCES

[1] ADLER, R. J. (1990). *An Introduction to Continuity, Extrema, and Related Topics for General Gaussian Processes.* IMS, Hayward, CA. MR1088478
[2] CHOI, T. and SCHERVISH, M. J. (2004). Posterior consistency in nonparametric regression problem under Gaussian process prior. Preprint.
[3] CHOUDHURI, N., GHOSAL, S. and ROY, A. (2004). Bayesian estimation of the spectral density of a time series. *J. Amer. Statist. Assoc.* **99** 1050–1059. MR2109494
[4] CHOUDHURI, N., GHOSAL, S. and ROY, A. (2007). Nonparametric binary regression using a Gaussian process prior. *Stat. Methodol.* **4**. To appear. Available at www4.stat.ncsu.edu/~sghosal/papers/binary_method.pdf.
[5] CRAMÉR, H. and LEADBETTER, M. R. (1967). *Stationary and Related Stochastic Processes. Sample Function Properties and Their Applications.* Wiley, New York. MR0217860
[6] GHOSAL, S., GHOSH, J. K. and RAMAMOORTHI, R. V. (1999). Posterior consistency of Dirichlet mixtures in density estimation. *Ann. Statist.* **27** 143–158. MR1701105
[7] GHOSAL, S., GHOSH, J. K. and VAN DER VAART, A. W. (2000). Convergence rates of posterior distributions. *Ann. Statist.* **28** 500–531. MR1790007
[8] HOEFFDING, W. (1963). Probability inequalities for sums of bounded random variables. *J. Amer. Statist. Assoc.* **58** 13–30. MR0144363
[9] LENK, P. J. (1988). The logistic normal distribution for Bayesian, nonparametric, predictive densities. *J. Amer. Statist. Assoc.* **83** 509–516. MR0971380
[10] LEONARD, T. (1978). Density estimation, stochastic processes and prior information (with discussion). *J. Roy. Statist. Soc. Ser. B* **40** 113–146. MR0517434
[11] TOKDAR, S. T. and GHOSH, J. K. (2007). Posterior consistency of Gaussian process priors in density estimation. *J. Statist. Plann. Inference* **137** 34–42.
[12] VAN DER VAART, A. W. and WELLNER, J. A. (1996). *Weak Convergence and Empirical Processes. With Applications to Statistics.* Springer, New York. MR1385671
[13] WALKER, S. G. (2004). New approaches to Bayesian consistency. *Ann. Statist.* **32** 2028–2043. MR2102501

DEPARTMENT OF STATISTICS
NORTH CAROLINA STATE UNIVERSITY
2501 FOUNDERS DRIVE
RALEIGH, NORTH CAROLINA 27695-8203
USA
E-MAIL: ghosal@stat.ncsu.edu

DEPARTMENT OF MATHEMATICS AND STATISTICS
UNIVERSITY OF MARYLAND, BALTIMORE COUNTY
1000 HILLTOP CIRCLE
BALTIMORE, MARYLAND 21250
USA
E-MAIL: anindya@math.umbc.edu